\documentclass[11pt]{amsart}

\usepackage{amsmath, amsfonts, amssymb,amsthm}
\usepackage{amstext}
\usepackage{mathrsfs}
\usepackage{float,epsf,subfigure}
\usepackage[all,cmtip]{xy}
\usepackage{hyperref}
\usepackage{algorithm}
\usepackage{algorithmic}
\usepackage{mathtools}
\usepackage{float}
\usepackage{bm}

\makeatletter
\@namedef{subjclassname@2020}{\textup{2020} Mathematics Subject Classification}
\makeatother

\newtheorem{theorem}{Theorem}[section]

\newtheorem{def-thm}[theorem]{Definition-Theorem}
\newtheorem{lemma}[theorem]{Lemma}

\newtheorem{defi}[theorem]{Definition}
\newtheorem{example}[theorem]{Example}
\newtheorem{remark}{Remark}

\newtheorem*{prop}{Proposition}
\newtheorem*{theo}{Theorem}
\newtheorem*{con1}{Yang's conjecture}
\newtheorem*{con2}{Conjecture}
\newtheorem*{tha}{Theorem A}

\theoremstyle{plain}
\theoremstyle{definition}

\begin{document}
\title[]{Periodicity of Entire Solutions of $ff''=H^2$}
\author[ ] {D.H. Zhao, S.S. Yang*, C.C. Yang and L.W. Liao}

\address{School of Mathematics \\ Nanjing University \\ Nanjing, P. R. China}
\email{donghaizhao2021@163.com}
\address{School of Mathematical Sciences\\ Suzhou University of Science and Technology\\ Suzhou, P. R. China}
\email{ssyangmath@126.com}

\address{Department of Mathematics\\The Hong Kong University of Science and Technology\\ Clearwater Bay, Kowloon, Hong Kong}
\email{maccyanlacy@gmail.com}

\address{School of Mathematics \\ Nanjing University \\ Nanjing, P. R. China}
\email{maliao@nju.edu.cn}

\thanks{This research is supported by the National Natural Science Foundation of China (Grant No. 12171232, 12471077, 12601138).}
\thanks{*Corresponding author}
\date{}

\begin{abstract}
Let $f$ denote a transcendental entire function. The so called Yang's conjecture asserts that if for some integer $k$ ($\geq1$) the product $ff^{(k)}$ is periodic, then $f$ must be also periodic. In this paper, we prove that when $k=2$ for a given periodic transcendental entire function $H$, and $f$ is an entire solution of the equation
$$ff''=H^2,$$
 then \(f\) must be periodic.
\end{abstract}
\maketitle \thispagestyle{empty} \setcounter{page}{1}

\vskip\baselineskip
\setlength\arraycolsep{2pt}
\medskip

\section{Introduction and main results}
First, we recall the following result, in which Laguerre's theorem (on the interlacing property of real zeros and critical points) was used:
\begin{tha} {\cite{Titchmarsh}}
The differential equation
\begin{equation}\label{1.1}
f(z)f''(z)=-\sin^2 z
\end{equation}
has no real solution, other than $f(z)=\pm\sin z$, which is an integral(entire) function of finite order.
\end{tha}
It is noted that
$$ff''=-\sin^2 z=(i\sin z)^2=H^2,$$
and $H$ is a periodic function.  Then a general question of this type was posed by C.C. Yang and also appeared as Conjecture 1.1 in \cite{WangHu2018}.

\begin{con1}
\emph{Let $f$ be a transcendental entire function, and let $k$ be a
positive integer. If $ff^{(k)}$ is periodic, then $f$ is periodic.}
\end{con1}

Wang and Hu \cite{WangHu2018} proved that $f$ is periodic if $(f^2)^{(k)}$ is periodic. Since $(f^2)'=2ff'$, their
result establishes Yang's conjecture for $k=1$.

By using Nevanlinna's theory, we have been able to prove the following result:
\begin{theo}
Let $H\not\equiv0$ be a periodic entire function.
If an entire function $f$ satisfies
$$ff''=H^2,$$
then $f$ is periodic.
\end{theo}
In Section 4, we apply the Theorem with $H(z)=i\sin z$ to give a new proof of Theorem~A.

\begin{remark}
For the Yang's conjecture, we believe that when $k>2$, by similarly but more complicated arguments, we can deduce the same conclusion.
\end{remark}

\par Now when $H$ is not periodic, then we have the following examples:
\begin{example}\label{ex1}
Let
$$f(z)=e^{\sqrt{2}z}\sinh^2 z,$$
and
$$H(z)=\sqrt{2}e^{\sqrt{2}z}\sinh z(\cosh z+\sqrt{2}\sinh z),$$
where $\rho(H)=1$ and $H$ is not periodic. Then
$$ff''=H^2.$$

\end{example}

\begin{example}\label{ex2}
Let $m\geq 1$ be an integer and $f$ be a nontrivial entire solution of
$$f''+z^{2m}f=0.$$
Then
$$ff''=-z^{2m}f^2=H^2,$$
where $H(z)=iz^mf(z)$. Consequently, $\rho(H)>1$ and $H$ is not periodic.
\end{example}
We also consider the more general equation
\begin{equation}\label{1.0}
    ff^{(k)}=H^n,
\end{equation}
where $k$ and $n$ are positive integers. The next example gives an entire solution for $k=1$ and every integer $n\ge3$.
\begin{example}\label{ex3}
Let $n\geq 3$ and set
$$f(z)=\exp\left(\frac{z^{n+1}}{n+1}\right),\quad H(z)=z\exp\left(\frac{2z^{n+1}}{n(n+1)}\right).$$
Then $H$ is not periodic and $f,H$ satisfy the equation $ff'=H^n$.

\end{example}
\par In each of the above examples, $H$ is a transcendental entire function with a zero. Examples \ref{ex1} and \ref{ex2} have $k=n=2$, while Example \ref{ex3} has $k=1$ and $n\ge3$. However, we have not found such an example with both $k\geq2$ and $n\ge3$. So, we may then raise the following conjecture:
\begin{con2}
Let $H$ be a given transcendental entire function. Then the differential equation \eqref{1.0} has no entire solutions if $k\geq2,n\geq3$ and $H$ has a zero.
\end{con2}
The following proposition supports our conjecture when $H$ has exactly one zero.
\begin{prop}
Let $H$ be a transcendental entire function having exactly one zero, and let $k\geq 2$ and $n\geq 1$.
Then the differential equation \eqref{1.0} has no entire solutions.
\end{prop}

\section{Some lemmas}

\begin{lemma}\label{L.Laine}{\rm\cite[p.35]{Laine}}
Let $f$ be a transcendental meromorphic function and $k\geq 1$ be an integer. Then
$$m\big(r,\frac{f^{(k)}}{f}\big)=O\big(\log T(r,f)+\log r\big)$$ possibly outside a set E of r with finite linear measure,
and if $f$ is of finite order of growth, then
$$m\big(r,\frac{f^{(k)}}{f}\big)=O(\log r).$$
\end{lemma}

\begin{lemma}\label{M-T}{\rm\cite{Laine}}
Let $g$ be an entire function and assume that $0<r<R<\infty$ and that the maximum modulus $M(r,g)=\max_{|z|=r}|g|$ satisfies $M(r,g)\geq1$. Then
$$T(r,g)\leq\log M(r,g)\leq\frac{R+r}{R-r}T(R,g).$$
\end{lemma}

\begin{lemma}\label{hayman1}{\rm\cite{Hayman1959}}
The only functions $f(z)$ meromorphic in the plane, such that $f$, $f''$ have only a finite number of zeros and poles are of the form
$$f(z)=\frac{P_1(z)}{P_2(z)}e^{P_3(z)}$$
where $P_1(z)$, $P_2(z)$, $P_3(z)$ are polynomials. Of these the only functions for which $f$ and $f''$ have no zeros are $e^{az+b}$ and $(az + b)^{-n}$, where $n$ is a positive integer.
\end{lemma}

\begin{defi}
Let $\rho>0$ and $0\leq\tau<\infty$. An entire function $f$ is said to be of \emph{growth $(\rho,\tau)$} if its order does not exceed $\rho$ and, when its order equals $\rho$, its type does not exceed $\tau$. In particular, an entire function of growth $(1,\tau)$ is said to be of \emph{exponential type $\tau$}.
\end{defi}

\begin{lemma}\label{Dry}{\rm\cite{Dry1}}
Let $f$ be a $2\pi$-$periodic$ entire function of exponential type $\tau$. Then, $f$ has the form
$$f(z):=\sum^{n}_{\nu=-n}c_{\nu}e^{i\nu z} \quad (n=[\tau]).$$
\end{lemma}

\begin{lemma}\label{hayman.con}{\rm\cite{Frank99}}
Suppose that $f$ is meromorphic in the plane and that $f$ and $f^{(k)}$ have only finitely many zeros, for some $k \geq2$. Then we have $f(z)=R(z)e^{P(z)}$, with $R$ a rational function and $P$ a polynomial. In particular, $f$ has finite order and finitely many poles.
\end{lemma}
\section{Proof of Theorem}
\begin{proof}
Let $\omega$ be a period of $H$. Then $P=H^2\not\equiv0$ is entire and $\omega$-periodic, and all its zeros are multiple.\par
If $P$ has no zero, Lemma \ref{hayman1} gives $f(z)=Ce^{az}$ with $aC\ne0$. Now we can suppose that $P$ has at least one zero. \par Let $\{u_j\}$ be a family of entire functions such that
$$u_j(z)=f(z+j\omega),\qquad  j\in\mathbb{Z},$$
it follows that $u_ju_j''=P$. Next we can assume that for any $i\not=j$, $u_i\not\equiv u_j$. \\
Define
\begin{equation}\label{3.0}
h_{ji}=u_j/u_i \quad and \quad w_i=u_i'/u_i \quad for \quad any \quad i\not=j\in \mathbb{Z}.
\end{equation}
 Hence we obtain
\begin{equation}\label{3.1}
\left(w_i\right)'=\frac{u_i''}{u_i}-\left(w_i\right)^2,
\end{equation}
and
\begin{equation}\label{3.2}
u_j''=h_{ji}''u_i+2h_{ji}'u_i'+h_{ji}u_i''.
\end{equation}
From $u_iu_i''=u_ju_j''$, it follows that
\begin{equation}\label{3.3}
\frac{u_i''}{u_i}=\frac{h_{ji}''h_{ji}}{1-h_{ji}^2}+2\frac{h_{ji}'h_{ji}}{1-h_{ji}^2}\frac{u_i'}{u_i}.
\end{equation}
Substituting \eqref{3.1} into \eqref{3.3}, we obtain
\begin{equation}\label{3.4}
\left(w_i\right)^2=\alpha_1+\alpha_2w_i-\left(w_i\right)' ,
\end{equation}
where
$$\alpha_1=\frac{h_{ji}''h_{ji}}{1-h_{ji}^2}=\frac{h_{ji}'h_{ji}}{1-h_{ji}^2}\frac{h_{ji}''} {h_{ji}'},$$
 and
$$\alpha_2=2\frac{h_{ji}'h_{ji}}{1-h_{ji}^2}=-\left(\frac{h_{ji}'}{h_{ji}-1}+\frac{h_{ji}'}{h_{ji}+1}\right).$$
Hence Lemma \ref{L.Laine} yields that
\begin{equation}\label{3.5}
m(r,\alpha_k)=O\big(\log T(r,h_{ji})+\log r\big), \quad k=1,2.
\end{equation}
Since $u_iu_i''=u_ju_j''$ we obtain
\begin{equation}\label{3.6}
T(r,h_{ji}^2)=T\left(r,\left(\frac{u_j}{u_i}\right)^2\right)=T\left(r,\frac{u_i''}{u_i}\frac{u_j}{u_j''}\right)\leq T\left(r,\frac{u_i''}{u_i}\right)+T\left(r,\frac{u_j''}{u_j}\right)+O(1).
\end{equation}
From \eqref{3.1} and \eqref{3.6} it follows that
$$2T(r,h_{ji})\leq T\left(r,\frac{u_i''}{u_i}\right)+T\left(r,\frac{u_j''}{u_j}\right)=O\left(T\left(r,w_i\right)+T\left(r,w_j\right)+\log r\right).$$
Therefore \eqref{3.5} yields
$$m(r,\alpha_k)=O\left(\log T\left(r,w_i\right)+\log T\left(r,w_j\right)+\log r\right), \quad k=1,2.$$
By \eqref{3.4}, we obtain
\begin{equation}\label{3.7}
\begin{aligned}
m\left(r,w_i\right)&=\frac{1}{2\pi}\int_0^{2\pi}\log^{+}\left|w_i(re^{i\theta})\right|\,d\theta\\
&=\frac{1}{2\pi}\int_{|w_i|\leq1}\log^{+}\left|w_i(re^{i\theta})\right|\,d\theta+\frac{1}{2\pi}\int_{|w_i|>1}\log^{+}\left|w_i(re^{i\theta})\right|\,d\theta\\
&=\frac{1}{2\pi}\int_{|w_i|>1}\log^{+}\left|\frac{\alpha_1}{w_i}+\alpha_2-\frac{w_i'}{w_i}\right|\,d\theta\\
&\leq m(r,\alpha_1)+m(r,\alpha_2)+m\left(r,\frac{w_i'}{w_i}\right)\\
&=O\left(\log T\left(r,w_i\right)+\log T\left(r,w_j\right)+\log r\right).
\end{aligned}
\end{equation}
The same conclusion holds after interchanging $i$ and $j$, hence \eqref{3.7} gives
\begin{equation}\label{3.8}
m\left(r,w_i\right)+m\left(r,w_j\right)=O\left(\log T\left(r,w_i\right)+\log T\left(r,w_j\right)+\log r\right).
\end{equation}
By \eqref{3.0}, every pole of $w_i$ or $w_j$ is simple and is a zero of $P$. Hence
$$N\left(r,w_i\right)+N\left(r,w_j\right)\leq2\overline{N}\left(r,\frac1P\right).$$
Combining this with \eqref{3.8}, we have
$$T\left(r,w_i\right)+T\left(r,w_j\right)\leq2\overline{N}\left(r,\frac1P\right)+O\left(\log T\left(r,w_i\right)+\log T\left(r,w_j\right)+\log r\right),$$
it follows that
\begin{equation}\label{3.9}
T\left(r,w_i\right)+T\left(r,w_j\right)\leq O\left(\overline{N}\left(r,\frac1P\right)+\log r\right),  \quad r\to\infty.
\end{equation}
Substituting \eqref{3.9} into \eqref{3.8} gives
\begin{equation}\label{3.10}
m\left(r,w_i\right)+m\left(r,w_j\right)=O\left(\log \overline{N}\left(r,\frac1P\right)+\log r\right).
\end{equation}
And from \eqref{3.1}, it follows that
\begin{equation}\label{3.11}
m\left(r,\frac{u_i''}{u_i}\right)=m(r,w_i'+w_i^2)\leq m\left(r,\frac{w_i'}{w_i}\right)+3m(r,w_i).
\end{equation}
Since Lemma \ref{L.Laine} and \eqref{3.9} give
$$m\left(r,\frac{w_i'}{w_i}\right)=O(\log T(r,w_i)+\log r)=O\left(\log \overline{N}\left(r,\frac1P\right)+\log r\right),$$
it follows from \eqref{3.10} and \eqref{3.11} that
\begin{equation}\label{3.12}
m\left(r,\frac{u_i''}{u_i}\right)=O\left(\log \overline{N}\left(r,\frac1P\right)+\log r\right).
\end{equation}
Here these estimates hold for all $i\not=j\in\mathbb{Z}$.\par
Next we prove that for any $i\not=j$, $u_i$ and $u_j$ have no common multiple zero. If not, we assume that $u_i$ and $u_j$ have zeros of multiplicity at least 2 at $z_0$. Hence both $u_i+u_j$ and $u_i-u_j$ would have zeros of multiplicity at least 2 at $z_0$. From $u_iu_i''=u_ju_j''$, it follows that
\begin{equation}\label{3.13}
(u_i+u_j)(u_i-u_j)''+(u_i-u_j)(u_i+u_j)''=2(u_iu_i''-u_ju_j'')=0.
\end{equation}
 We can express $u_i+u_j$ and $u_i-u_j$ as Taylor series expansions around $z_0$:
 \begin{equation}\label{3.14}
\begin{aligned}
u_i+u_j&=A(z-z_0)^\alpha+O((z-z_0)^{\alpha+1}),\\
u_i-u_j&=B(z-z_0)^\beta+O((z-z_0)^{\beta+1}),
\end{aligned}
\end{equation}
 where $A,B\not=0$ and $\alpha,\beta\geq2$. Based on \eqref{3.14}, we can derive the following expressions for $(u_i+u_j)''$ and $(u_i-u_j)''$:
  \begin{equation}\label{3.15}
\begin{aligned}
(u_i+u_j)''&=A\alpha(\alpha-1)(z-z_0)^{\alpha-2}+O((z-z_0)^{\alpha-1}),\\
(u_i-u_j)''&=B\beta(\beta-1)(z-z_0)^{\beta-2}+O((z-z_0)^{\beta-1}).
\end{aligned}
\end{equation}
By substituting expressions \eqref{3.14}-\eqref{3.15} into equation \eqref{3.13}, we obtain
$$AB(\alpha(\alpha-1)+\beta(\beta-1))(z-z_0)^{\alpha+\beta-2}+O((z-z_0)^{\alpha+\beta-1})=0,$$
which yields $AB(\alpha(\alpha-1)+\beta(\beta-1))=0$. Since $A,B\not=0$ and $\alpha,\beta\geq2$, a contradiction. Thus at most one of $\{u_j\}$ has a multiple zero at a fixed point $z_0$.

Now we consider an auxiliary function
\begin{equation}\label{3.16}
\Phi=\left(\prod_{j=0}^7\frac{u_j''}{u_j}\right)^6\prod_{0\le i<j\leq7}\left(w_i-w_j\right),
\end{equation}
where $\Phi\not\equiv0$, otherwise $\Phi\equiv0$ yields $w_i-w_j\equiv0$. It follows that
$$\left(\frac{u_i}{u_j}\right)'=\frac{u_i'u_j-u_iu_j'}{u_j^2}=\frac{u_i}{u_j}(w_i-w_j)\equiv0,$$
which gives that $u_i/u_j$ would be a nonzero constant. By $u_iu_i''=u_ju_j''$, we obtain $u_i=u_j$, a contradiction.\par\
From \eqref{3.10}, \eqref{3.12} and \eqref{3.16}, it follows that
\begin{equation}\label{3.17}
\begin{aligned}
m(r,\Phi)&\leq O\left(\sum_{j=0}^7m\left(r,\frac{u_j''}{u_j}\right)+\sum_{j=0}^7m(r,w_j)\right)\\
&=O\left(\log \overline{N}\left(r,\frac1P\right)+\log r\right).
\end{aligned}
\end{equation}
\par We next claim that $\Phi$ is entire and all zeros of $P$ are zeros of $\Phi$. Let $z_0$ be a zero of $P$ with multiplicity $m_p\geq2$. And for each $j\in\{0,\dots,7\}$ we consider the following cases:

\smallskip\noindent\emph{Case 1: $u_j(z_0)\ne0$.}

Since
$$ \frac{u_j''}{u_j}=\frac{P}{u_j^2},$$
we obtain $z_0$ is a zero of $ {u_j''}/{u_j}$ with multiplicity $m_p$. And $w_j=u_j'/u_j$ is regular at $z_0$.

\smallskip\noindent\emph{Case 2: $z_0$ is a simple zero of $u_j$.}

From $u_j u_j''=P$, it follows that $u_j''/u_j$ is regular and has a zero of multiplicity $m_p-2\geq0$ at $z_0$. We can express $u_j$ as a Taylor expansion around $z_0$:
\begin{equation}\label{3.18}
u_j(z)=c_1(z-z_0)+c_2(z-z_0)^2+O\left((z-z_0)^3\right),
\end{equation}
where $c_1\not=0$. Based on \eqref{3.18}, we can derive the following expressions for the derivatives of $u_j$:
\begin{equation}\label{3.19}
\begin{aligned}
u_j'(z)&=c_1+2c_2(z-z_0)+O\left((z-z_0)^2\right),\\
u_j''(z)&=2c_2+O(z-z_0).
\end{aligned}
\end{equation}
Since $m_p\geq2$, $u_j u_j''=P$ yields $u_j''$ has a zero at $z_0$ of multiplicity $m_p-1\geq1$. It follows from \eqref{3.19} that $c_2=0$.
Hence, \eqref{3.18}-\eqref{3.19} can reduce to
\begin{equation}\label{3.20}
\begin{aligned}
u_j(z)&=c_1(z-z_0)+O\left((z-z_0)^3\right),\\
u_j'(z)&=c_1+O\left((z-z_0)^2\right),\\
u_j''(z)&=O(z-z_0).
\end{aligned}
\end{equation}
Based on \eqref{3.20}, we can derive the following expression for the logarithmic derivative of $u_j$:
\begin{equation}\label{3.21}
w_j=\frac{u_j'}{u_j}=\frac{1}{z-z_0}+O(z-z_0).
\end{equation}
\smallskip\noindent\emph{Case 3: $z_0$ is a zero of $u_j$ with multiplicity $m\geq2$.}

We can express $u_j$ as a Taylor expansion around $z_0$:
\begin{equation}\label{3.22}
u_j(z)=c_m(z-z_0)^m+O\left((z-z_0)^{m+1}\right),
\end{equation}
where $c_m\ne0$. Based on \eqref{3.22}, we can derive the following expressions for the derivative of $u_j$ and its logarithmic derivative:
\begin{equation}\label{3.23}
\begin{aligned}
u_j'(z)&=mc_m(z-z_0)^{m-1}+O\left((z-z_0)^{m}\right),\\
w_j&=\frac{u_j'}{u_j}=\frac{m}{z-z_0}+O(1).
\end{aligned}
\end{equation}
Since $z_0$ is a zero of $u_j$ with multiplicity $m\geq2$, we obtain $u_j''/u_j$ has a pole at $z_0$ of multiplicity 2.
\par The definition \eqref{3.16} yields that if $\Phi$ has a pole at $z_1$ this point must be a zero of $u_j$ for some $j\in\{0,\dots,7\}$. And from $u_ju_j''=P(z)$, it follows that $z_1$ must be a zero of $P$. Let $\lambda$ be the number of $u_j,j\in\{0,\dots,7\}$ which satisfies $Case$ 1, and let $\mu$ be the number of $u_j,j\in\{0,\dots,7\}$ which does not satisfy $Case$ 3. Hence the zero multiplicity of
$$\prod_{j=0}^7\frac{u_j''}{u_j}$$
at $z_0$ is
$$m_1=\lambda m_p+(\mu-\lambda)(m_p-2)-2(8-\mu)=\mu m_p+2\lambda-16.$$
\par Next we estimate the zero multiplicity of
$$\prod_{0\le i<j\leq7}\left(w_i-w_j\right).$$
If one of $u_i,u_j$ satisfies $Case$ 3, it follows from \eqref{3.21} and \eqref{3.23} that $w_i-w_j$ has a simple pole at $z_0$. \\
If all of $u_i,u_j$ satisfy $Case$ 2, it follows from \eqref{3.21} that
$$w_i-w_j=O(z-z_0),$$
which means that $w_i-w_j$ has a zero at $z_0$ with multiplicity at least 1.\\
If one of $u_i,u_j$ satisfies $Case$ 2 and another satisfies $Case$ 1, it follows from \eqref{3.21} that $w_i-w_j$ has a simple pole at $z_0$.\\
If all of $u_i,u_j$ satisfy $Case$ 1, $w_i-w_j$ is regular at $z_0$.\\
Hence we obtain that the zero multiplicity of
$$\prod_{0\le i<j\leq7}\left(w_i-w_j\right)$$
 at $z_0$ is
 $$m_2\geq-\mu(8-\mu)+\binom{\mu-\lambda}{2}-\lambda(\mu-\lambda).$$
From at most one of $\{u_j\}$ has a multiple zero at a fixed point $z_0$, it follows that $7\leq\mu\leq8$. Hence we obtain
\begin{equation}\label{3.24}
\begin{aligned}
6m_1+m_2&\geq6(8m_p+2\lambda-16)+\binom{8-\lambda}{2}-\lambda(8-\lambda)\\
&\geq 12\lambda+\binom{8-\lambda}{2}-\lambda(8-\lambda)\\
&=\frac{3\lambda^2-7\lambda+56}{2}>0,
\end{aligned}
\end{equation}
 as $\mu=8$, and
 \begin{equation}\label{3.25}
\begin{aligned}
6m_1+m_2&\geq6(7m_p+2\lambda-16)-7+\binom{7-\lambda}{2}-\lambda(7-\lambda)\\
&\geq 6(2\lambda-2)-7+\binom{7-\lambda}{2}-\lambda(7-\lambda)\\
&=\frac{3\lambda^2-3\lambda+4}{2}>0,
\end{aligned}
\end{equation}
 as $\mu=7$.
 \par Hence \eqref{3.24}-\eqref{3.25} gives that every zero of $P$ is a zero of $\Phi$, and $\Phi$ is entire. It follows from \eqref{3.17} that
$$\begin{aligned}
\overline{N}\left(r,\frac1P\right)
&\leq N\left(r,\frac1\Phi\right)\\
&\leq T(r,\Phi)+O(1)\\
&=m(r,\Phi)+O(1)\\
&=O\left(\log\overline{N}\left(r,\frac1P\right)+\log r\right).
\end{aligned}$$
Therefore we obtain
$$\overline{N}\left(r,\frac1P\right)=O(\log r).$$
On the other hand, since $P$ has at least one zero and $P$ is $\omega$-periodic, we can obtain there is a constant $c$ such that
$$\overline{N}\left(r,\frac1P\right)>c r,$$
as $r\to\infty$, a contradiction. Hence there are $i,j$ such that $i\not=j$ and $u_i\equiv u_j$, which yields $f$ is periodic.
\end{proof}

\section{A New Proof of Theorem A}
\begin{proof}
Let $f$ be an entire solution of \eqref{1.1} and take $H=i\sin z$. Then from our theorem, we conclude that $f$ is periodic. If $\omega$ is a non-zero period of $f$, then it is a period of $\sin^2 z$. Hence $\omega\in\pi\mathbb Z$ and $f$ has a period $2N\pi$ for some $N\in\mathbb{Z_+}$.
\par Set $G=f''/f$. Since every zero of $f$ is a zero of $\sin z$ and every pole of $G$ has multiplicity at most 2, we obtain
$$N(r,G)\leq 2N\left(r,\frac1{\sin z}\right)=O(r).$$
Since Lemma \ref{L.Laine} yields
$$m(r,G)=O(\log T(r,f)+\log r).$$
It follows from $f^2=-\sin^2 z/G$ that
$$2T(r,f)\leq 2T(r,\sin z)+T(r,G)+O(1)=O(r)+O(\log T(r,f)),$$
which means that
$$T(r,f)=O(r).$$
Hence Lemma \ref{M-T} gives $\log M(r,f)=O(r)$, so $f$ is of
finite exponential type.

Put $F(z)=f(Nz)$. Then Lemma \ref{Dry} yields
 \begin{equation}\label{1.2}
 F(z)=\sum_{\nu=p}^{q}c_\nu e^{i\nu z},
 \end{equation}
where $c_p,c_q\not= 0$.
From \eqref{1.1} it follows that
 \begin{equation}\label{1.3}
 FF''=-N^2\sin^2(Nz)=\frac{N^2}{4}\left(e^{2iNz}-2+e^{-2iNz}\right).
 \end{equation}
Comparing the exponential terms with the smallest and largest indices in \eqref{1.3} gives $p=-N$ and $q=N$. Let $-N<j<N$ be the largest integer for which $c_j\ne0$. From \eqref{1.2}, it follows that $FF''$ has an exponential term $e^{i(N+j)z}$ and its coefficient is $-(N^2+j^2)c_Nc_j\not=0$. Since $0< N+j<2N$, equation \eqref{1.3} forces $N+j=0$. Hence
$$F(z)=c_Ne^{iNz}+c_{-N}e^{-iNz},$$
and
$$F''(z)=-N^2c_Ne^{iNz}-N^2c_{-N}e^{-iNz},$$
which means that $F''=-N^2F$. From $F(z)=f(Nz)$, it follows that $f''=-f$. Therefore \eqref{1.1} gives $f^2=\sin^2z$. This proves Theorem A.
\end{proof}
\section{Proof of Proposition}
\begin{proof}
Suppose that $f$ is an entire solution of equation \eqref{1.0}.
\par We can assume that the only zero of $H$ is 0. It follows that every zero of $f$ or $f^{(k)}$ is 0. Since $f$ and $f^{(k)}$ have only finitely many zeros, then Lemma \ref{hayman.con} gives
\begin{equation}\label{5.1}
    f(z)=z^m e^{Q(z)},
\end{equation}
where $m\geq0$ and $Q$ is nonconstant. Let $q=\deg Q\geq1$. Then for any $1\leq t\leq k$, we obtain
\begin{equation}\label{5.3}
f^{(t)}=P_te^{Q(z)}
\end{equation}
where $P_0=z^m$, $P_t=P_{t-1}Q'+P_{t-1}'$ and $\deg P_t=m+t(q-1)$.
\par From every zero of $f^{(k)}$ is 0 it follows that
\begin{equation}\label{5.2}
    f^{(k)}(z)=bz^s e^{Q(z)},
\end{equation}
where $b$ is a constant and $s=m+k(q-1)$.
\par We next claim that $m=s=0$. From \eqref{5.3} we get
\begin{equation}\label{5.4}
P_k=z^mQ'^k+mk(z^{m-1})Q'^{k-1}+\frac{k(k-1)}{2}z^mQ'^{k-2}Q''+O(z^{s-q-1}).
\end{equation}
If $q=1$, then $s=m$ and \eqref{5.4} gives
$$P_k=z^mQ'^k+mk(z^{m-1})Q'^{k-1}+O(z^{m-2}),$$
it follows from $P_k=bz^s$ that $mkQ'^{k-1}=0$. Since $k,Q'$ are non-zero, we obtain $m=s=0$.\\
Next we consider the case $q\geq2$. Let $a\ne0$ be the leading coefficient of $Q'$. Then comparing the leading coefficients in $P_k=bz^s$ gives $b=a^k$. Since $q,k\geq2$, we have
$$ s-q=m+(k-1)(q-1)-1\ge0$$
and
$$mk+\frac{k(k-1)}2(q-1)>0.$$
Set $R=P_k-z^mQ'^k$, then \eqref{5.4} yields
$$R=a^{k-1}\left(mk+\frac{k(k-1)}2(q-1)\right)z^{s-q}+O(z^{s-q-1}),$$
where $R\not\equiv0$ and  $\deg R=s-q$.
On the other hand, since $P_k=bz^s=a^kz^s$, it follows that
$$\begin{aligned}
R&=P_k-z^m(Q')^k\\
&=z^m\left((az^{q-1})^k-(Q')^k\right)\\
&=z^m(az^{q-1}-Q')\sum_{j=0}^{k-1}(az^{q-1})^{k-1-j}(Q')^j.
\end{aligned}
$$
Hence $R$ has only one highest degree term
$$z^m(az^{q-1}-Q') ka^{k-1}z^{(k-1)(q-1)},$$
which means that
$$\begin{aligned}
    \deg R
    &=m+\deg(az^{q-1}-Q')+(k-1)(q-1)\\
    &\ge m+(k-1)(q-1)\\
    &=s-q+1,
\end{aligned}$$
a contradiction. Hence we obtain $m=s=0$. By substituting \eqref{5.1} and \eqref{5.2} into \eqref{1.0}, we obtain
$$ff^{(k)}=be^{2Q}=H^n.$$
Since $H$ has exactly one distinct zero, a contradiction. This completes the proof.
\end{proof}

\vskip\baselineskip

\end{document}